\input amstex 
\documentstyle{amsppt}
\input bull-ppt
\keyedby{bull348/lbd}
\define\vare{\varepsilon}

\topmatter
\cvol{28}
\cvolyear{1993}
\cmonth{January}
\cyear{1993}
\cvolno{1}
\cpgs{99-103}
\title M\"obius invariance of knot energy\endtitle
\author Steve Bryson, Michael H. Freedman, Zheng-Xu He, 
and Zhenghan Wang
\endauthor
\shortauthor{Steve Bryson, M. H. Freedman, Z.-X. He, and 
Zhenghan Wang}
\shorttitle{M\"obius invariance of knot energy}
\address (Steve Bryson)
Computer Sciences Corporation, NASA Ames Research Center,
Moffett Field, California 94035\endaddress
\ml bryson\@nas.nasa.gov\endml
\address (Michael H. Freedman, Zhenghan Wang)
Department of Mathematics, University of California at San 
Diego, La
Jolla, California 92093-0112\endaddress
\ml M. H. Freedman mfreedman\@ucsd.edu\endml
\address (Zheng-Xu He) Department of Mathematics, 
Princeton University,
Princeton, New Jersey
08544-1000\endaddress
\curraddr Department of Mathematics, Cornell University, Ithaca, New 
York 14853\endcurraddr
\date April 14, 1992\enddate
\thanks The first author thanks the San Diego 
Supercomputer Center for use of
their facilities. The second and fourth authors were 
supported in part by NSF
grant DMS-8901412. The third author was supported in part 
by NSF grant
DMS-9006954\endthanks
\abstract A physically natural potential energy for simple 
closed curves in
$\bold R^3$ is shown to be invariant under M\"obius 
transformations. This leads
to the rapid resolution of several open problems: round 
circles are precisely
the absolute minima for energy; there is a minimum energy 
threshold below which
knotting cannot occur; minimizers within prime knot types 
exist and are
regular. Finally, the number of knot types with energy 
less than any constant
$M$ is estimated.\endabstract
\endtopmatter

\document
Consider a rectifiable curve $\gamma(u)$ in the Euclidean 
3-space $\bold R^3$,
where $u$ belongs to $\bold R^1$ or $S^1$. Define its 
energy by
$$E(\gamma)=\iint\left\{\frac{1}{|\gamma(u)-\gamma(v)|^2}-%
\frac{1}
{D(\gamma(u),\gamma(v))^2}\right\}|\dot{\gamma}(u)||\dot{%
\gamma}(v)|\,du\,dv,$$
where $D(\gamma(u),\gamma(v))$ is the shortest arc 
distance between$\gamma(u)$
and $\gamma(v)$ on the curve. The second term of the 
integrand is called a
regularization (see \cite{O1--O3, FH}). It is easy to see 
that $E(\gamma)$ is
independent of parametrization and is unchanged if 
$\gamma$ is changed by a
similarity of $\bold R^3$.
\par
Recall that the M\"obius transformations of the 3-sphere 
$=\bold R^3\cup\infty$
are the ten-dimensional group of angle-preserving 
diffeomorphisms generated by
inversion in 2-spheres.
\par
The central fact of this announcement is:
\thm{M\"obius Invariant Property} Let $\gamma$ be a closed 
curve in $\bold
R^3$. If $T$ is a M\"obius transformation of $\bold 
R^3\cup\infty$ and
$T(\gamma)$ is contained in $\bold R^3$, then 
$E(T(\gamma))=E(\gamma)$. If
$T(\gamma)$ passes through $\infty$, the integral satisfies
$E(T(\gamma))=E(\gamma)-4$.
\ethm
\par
This simple fact (proved below), combined with earlier 
results proved in
\cite{FH}, allows the rapid resolution of several open 
problems.
\thm{Theorem A} Among all rectifiable loops $\gamma\colon\ 
S^1\to\bold R^3$,
round circles have the least energy {\rm(}E {\rm(round 
circle)} $=4)$ and any
$\gamma$ of least energy parameterizes a round circle.
\ethm
\thm{Theorem B} If $K$ is a smooth prime {\rm(}not a 
connected sum{\rm)} knot,
then there exists a simple closed rectifiable $\gamma_K$ 
of knot type $K$ with
$E(\gamma_K)\le E(\gamma)$ for all rectifiable loops 
$\gamma$ which are
topologically ambient isotopic to $K$.
\ethm
\thm{Theorem C} Any minimizer $\gamma_K$, as above, will 
enjoy some
regularity. With an arc length parametrization, $\gamma_K$ 
will be in
$C^{1,1}$.
\ethm
\par
Several results of \cite{FH} can be improved quantitatively.
\thm{Theorem D} If $\gamma$ is topologically tame, let 
$c([\gamma])$ denote the
{\rm(}topological{\rm)} crossing number of the knot type. 
We have
$$2\pi c([\gamma])+4\le E(\gamma).$$
{\rm(}It was proved in \cite{FH} that finite energy 
implies tame.{\rm)}
\ethm
\par
Since an essential knot must have three or more crossings, 
we obtain the
following
\thm{Corollary} Any rectifiable loop with energy less than 
$6\pi+4\approx
22.84954$ is unknotted.
\ethm
\par
Computer experiments of \cite{A} as reported in \cite{O3} 
and independently by
the first author yield an essential knot (a trefoil) with 
energy $\approx74$.
\par
It may be estimated \cite{S,T,W} that the number $K(n)$ of 
distinct knots of
at most $n$ crossings satisfies
$$2^n\le K(n)\le 2\cdot 24^n.$$
Hence the number of knot types with representatives below 
a given energy
threshold can also be bounded by an exponential.
\thm{Theorem E} The number $K_e(M)$ of isomorphism classes 
of knots which have
representatives of energy less than or equal to
$M$ is bounded by $2(24^{-4/2\pi})(24^{1/2\pi})^M
\approx
(0.264)(1.658)^M$. In particular, only finitely many
knot types occur below any finite energy threshold.
\ethm
\par
Note that there are competing candidates for the exponent 
$=-2$ in the
definition of $E$; for example, the Newtonian potential in 
$\bold R^3$ has
exponent $=-1$. When the exponent is strictly larger 
than$-3$, finite values
are obtained for smooth simple loops. Exponents smaller or 
equal to $-2$ yield
energies which blow up as a simple loop $\gamma$ begins to 
acquire a double
point, thus creating an infinite energy barrier to a 
change of topology. Such a
barrier would not exist for the Newtonian potential. We 
refer to \cite{O1--O3}
for detailed discussions. Similarity and M\"obius 
invariance are, of course,
special to the exponent $-2$.
\demo{Proof of Theorem {\rm A}} Let $T$ be a M\"obius 
transformation sending a point
of $\gamma$ to infinity. The energy $E(T(\gamma))\ge0$ 
with equality holding
iff $T(\gamma)$ is a straight line. Apply the M\"obius 
invariant property to
complete the proof.\qed\enddemo
\demo{Proof of Theorem {\rm B}} In \cite{FH} it is shown 
that for prime knot types
$K$ minimizers exist in the class of properly embedded 
rectifiable lines whose
completion in $\bold R^3\cup\infty$ represent $K$. 
According to the M\"obius
Invariance Property, such lines may be moved to a closed 
minimizer by any
M\"obius transformation $T$ which moves the completed line 
off infinity.\qed\enddemo
\demo{Sketch of Proof of Theorem {\rm C}} Let $\gamma_K$ 
be a closed minimizer in
knot type $K$. An inversion argument shows that, for 
sufficiently small
$\vare>0$, if $\gamma_K$ meets a closed ball of radius 
$\vare$, $B_\vare$, only
in its boundary $S_\vare$, then $\gamma_K\cap S_\vare$ 
consists of (at most)
one point. The idea is that if $\gamma_K\cap S_\vare$ is 
disconnected,
inverting an arc of $\gamma_K\backslash S_\vare$ into 
$B_\vare$ will lower
energy while preserving the knot type. Thus there is a 
continuous projection
from the $\vare$-neighborhood of $\gamma_K$ to $\gamma_K$ 
given by ``closest
point'' $\pi\colon\ \scr N_\vare(\gamma_K)\to\gamma_K$. We 
prove that the fibers
$\pi^{-1}(pt)$ are all geometric planar disks of radius 
$\vare$. The
disjointness of these ``normal'' fibers to distance 
$\vare$ is equivalent to the
existence of a continuously turning tangent to $\gamma_k$ 
whose generalized
derivative is in $L^\infty$.\qed\enddemo
\par
A detailed proof of Theorem C will appear elsewhere.
\demo{Proof of Theorem {\rm D}} Theorem 2.5 of \cite{FH} 
gives the inequality
$$c([\gamma])\le c(\gamma)\le E(\gamma)/2\pi$$
for proper rectifiable lines. According to the M\"obius 
Invariance Property,
the energy will increase by exactly 4 if a M\"obius 
transformation is used to
move the line off infinity and into closed position.\qed\enddemo
\demo{Proof of M\"obius Invariance Property} It is 
sufficient to consider how
$I$, an inversion in a sphere, transforms energy. Let $u$ 
be the arc length
parameter of a rectifiable closed curve $\gamma,u\in\bold 
R/l\bold Z$. Let
$$E_\vare(\gamma)=\iint_{|u-v|\ge\vare}\left(\frac{1}{|%
\gamma(u)-\gamma(v)|^2}-
\frac{1}{(D(\gamma(u),\gamma(v)))^2}\right)\,du\,dv\tag1$$
and
$$\aligned
E_\vare(I\circ\gamma)=&\iint_{|u-v|\ge\vare}
\left(\frac{1}{|I\circ\gamma(u)-I\circ\gamma(v)|^2}
-\frac{1}{(D(I\circ\gamma(u),I\circ\gamma(v)))^2}\right)\\
&\qquad\times\|I'(\gamma(u))\|\cdot\|I'(\gamma(v))\|\,du%
\,dv.\endaligned\tag2$$
\par
Clearly $E(\gamma)=\lim_{\vare\to0}E_\vare(\gamma)$ and
$E(I\circ\gamma)=\lim_{\vare\to0}E_\vare(I\circ\gamma)$.
\par
It is a short calculation (using the law of cosines) that 
the first terms
transform correctly, i.e.,
$$\frac{\|I'(\gamma(u))\|\cdot\|I'(\gamma(v))\|}{|I(%
\gamma(u))-I(\gamma(v))|^2}=
\frac{1}{|\gamma(u)-\gamma(v)|^2}.$$
\par
Since $u$ is arclength for $\gamma$, the regularization 
term of (1) is the
elementary integral
$$\int_{u=0}^l\left[2\int_{v=\vare}^{l/2}\frac{1}{v^2}\,dv%
\right]\,du=4-
\frac{2l}{\vare}.\tag3$$
\par
Let $s$ be an arclength parameter for $I\circ\gamma$.
Then $ds(u)/du=\|I'(\gamma(u))\|$ where 
$\|I'(\gamma(u))\|=f(u)$ denotes the
linear expansion factor of $I'$. Since $\gamma(u)$ is a 
lipschitz function and
$I'$ is smooth, $f(u)$ is lipschitz, hence, it has a 
generalized derivative
$f'(u)\in L^\infty$.
$$\aligned
\roman{regularization}\ (2)=&\int_{u\in\bold R/l\bold
Z}\left[\int_{|v-u|\ge\vare}\frac{|(I\circ\gamma)'(v)|\,dv}
{D(I\circ\gamma(u),I\circ\gamma(v))^2}\right]|(I\circ%
\gamma)'(u)|\,du\\
=&\int_{\bold R/l\bold Z}\left[\frac{4}{L}-\frac{1}{\vare_+
}-\frac{1}{\vare_-}
\right]\,ds,\endaligned\tag4$$
where $L=\roman{Length}(I(\gamma))$ and
$$\align
\vare_+&=\vare_+(u)=D((I\circ\gamma)(u),(I\circ\gamma)
(u+\vare))=s(u+\vare)-s(u)\\
&=\int_u^{u+\vare}f(w)\,dw
=f(u)\vare+\vare^2\int_0^1(1-t)f'(u+\vare t)\,dt\endalign$$
and
$$
\vare_-=\vare_-(u)=D((I\circ\gamma)(u-\vare),(I\circ%
\gamma)(u))
=f(u)\vare-\vare^2\int_0^1(1-t)f'(u-\vare t)\,dt.$$
Since $|f'(w)|$ is uniformly bounded, we have
$$\align
\frac{1}{\vare_+}=&\frac{1}{f(u)\vare}\left[\frac{1}{1+
(\vare/f(u))
\int_0^1(1-t)f'(u+\vare t)\,dt}\right]\\
=&\frac{1}{f(u)\vare}\left[1-\frac{\vare}{f(u)}%
\int_0^1(1-t)f'(u+\vare t)\,dt+
\scr O(\vare^2)\right]\\
=&\frac{1}{f(u)\vare}-\frac{1}{f(u)^2}\int_0^1(1-t)f'(u+
\vare t)\,dt+\scr
O(\vare).\endalign$$
Similarly,
$$\frac{1}{\vare_-}=\frac{1}{f(u)\vare}+
\frac{1}{f(u)^2}\int_0^1(1-t)f'(u-\vare
t)\,dt+\scr O(\vare).$$
Then by (4)
$$\aligned
\roman{regularization}\ (2)=&4-\int_{\bold R/l\bold 
Z}\frac{2}{\vare}\,du\\
&+\iint_{\bold R/l\bold 
Z\times[0,1]}\frac{(1-t)}{f(u)}[f'(u+\vare t)-f'(u-
\vare t)]\,du\,dt+\scr O(\vare)\\
=&4-\frac{2l}{\vare}+\scr O(\vare)+\scr 
O(\vare).\endaligned\tag5$$
\par
Comparing (3) and (5), we get
$$E_\vare(\gamma)-E_\vare(I\circ\gamma)=\scr O(\vare);$$
hence, $E(\gamma)=E(I\circ\gamma)$.
\par
For the second assertion, let $I$ send a point of $\gamma$ 
to infinity. In this
case $L=\infty$ and, thus, the constant term 4 in (5) 
disappears.\qed\enddemo
\heading Acknowledgment\endheading
\par
The authors wish to thank Adriano Garsia and Fred Hickling 
for useful
discussions.

\enddocument